\documentclass[12pt]{article}

\usepackage{setspace}
\usepackage{latexsym}
\usepackage{amsmath}
\usepackage{epsfig}
\usepackage{epic}
\usepackage{amssymb}

\newtheorem{lemma}{Lemma}[section]
\newtheorem{theorem}[lemma]{Theorem}

\newtheorem{corollary}[lemma]{Corollary}

\newcommand{\pf}{\noindent{\em Proof: }}
\newcommand{\epf}{\hfill\hbox{\rule{3pt}{6pt}}\\}

\newcommand{\cP}{\mbox{$\mathcal P$}}

\newcommand{\cC}{\mathcal{C}}
\newcommand{\exc}{{\rm exc}}
\newcommand{\Exc}{{\rm Exc}}

\newcommand\ring[1]{\mathaccent23{#1}}
\def\Er{\ring{E}}

\begin{document} 
\title{Slim Sets of Binary Trees}
\author{\bf Stefan Gr{\"u}newald
\\CAS-MPG Partner Institute for Computational Biology 
\\Shanghai Institutes for Biological Sciences 
\\320 Yue Yang Road 
\\Shanghai 200031, China 
\\and 
\\Max Planck Institute for Mathematics in the Sciences 
\\Inselstra{\ss}e 22 
\\04103 Leipzig, Germany
\\stefan@picb.ac.cn
\\Tel: +86 (21) 54920459
\\Fax: +86 (21) 54920451}



\date{\today}
\maketitle

\newpage
\begin{abstract}
A classical problem in phylogenetic tree analysis is to decide
whether there is a phylogenetic tree $T$ that contains all
information of a given collection $\cP$ of phylogenetic trees. 
If the answer is ``yes'' we say that $\cP$ is compatible and $T$
displays $\cP$. 
This decision problem is NP-complete even if all input trees are
quartets, that is binary trees with exactly four leaves. In this
paper, we prove a sufficient condition for a set of binary
phylogenetic trees to be compatible. That result is used to give a short and
self-contained proof of the known characterization of quartet sets of
minimal cardinality which are displayed by a unique phylogenetic tree. 
\end{abstract}


\section{Introduction}
Phylogenetic trees are used in computational biology to visualize the
evolutionary relationship between some taxa (e.g. species or
genes). For many methods in phylogenetic tree reconstruction, the required
running time grows at least exponentially with the
number of taxa. Hence, it is a natural and widely used approach to construct
a collection of phylogenetic trees with small taxa sets first and then a
supertree (see \cite{b04} for an overview), that is a phylogenetic tree
for all taxa that are contained in at least one of the small trees. 
However, it is NP-complete to decide whether there is a supertree that 
displays all input trees, even if each of them contains only four taxa
\cite{s92}.

\begin{figure}[b]
\center
\input{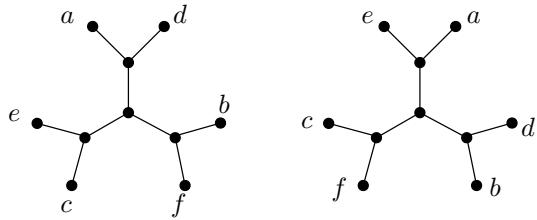}
\caption{Two phylogenetic trees.} 
\label{snow}
\end{figure}

Mathematically, a {\em phylogenetic tree} is a tree without
vertices of degree~2 and it is called {\em binary} if every interior
vertex has degree~3. The set of leaves or {\em taxa set} of a
phylogenetic tree $T$ is denoted by $L(T)$ and,
for a set $\cP$ of phylogenetic trees, we define $L(\cP)=
\displaystyle\bigcup_{T \in \cP}L(T)$. Let $T$ and $T'$ be two
phylogenetic trees with $L(T') \subseteq L(T)$ where $T'$ is
binary. We say that $T'$ is
{\em displayed} by $T$ if $T'$ is the phylogenetic tree obtained from the
smallest subtree of $T$ that contains $L(T')$ by suppressing all
vertices of degree~2. In that case $T'$ is called the {\em restriction
  of $T$ to $L(T')$} and is denoted by $T|_{L(T')}$. 
A set $\cP$ of binary phylogenetic trees is {\em
  compatible} if there is a phylogenetic tree $T$ that displays every tree
in $\cP$ and in that case we say that $T$ {\em displays} $\cP$. A {\em
  cherry} of a phylogenetic tree is a set of two leaves which are
adjacent to the same vertex of degree~3. A {\em quartet} is a binary
phylogenetic tree with exactly four leaves and the quartet with leaf
set $\{a,b,c,d\}$ and cherries $\{a,b\}$ and $\{c,d\}$ is denoted by
$ab|cd$. A set $\cP$ of binary phylogenetic trees is {\em definitive}
if there is exactly one
phylogenetic tree $T$ with taxa set $L(\cP)$ displaying $\cP$ and in
that case we say that $\cP$ {\em defines} $T$. Here, two phylogenetic
trees $T_1$ and $T_2$ are defined to be the same (or isomorphic) if
$L(T_1)=L(T_2)$ and there is a graph isomorphism from $T_1$ to $T_2$
that maps every leaf to itself. For example, both phylogenetic trees depicted 
in Figure \ref{snow} display the quartets $ad|cf$, $ae|bf$, and $bd|ce$ thus 
the set $\{ad|cf,ae|bf,bd|ce\}$ is not definitive. 

The main result of this paper is a sufficient condition for a set of
binary phylogenetic trees to be compatible. We need some more
definitions to state it. 
The set of interior edges of a phylogenetic tree $T$ is denoted by
$\Er(T)$. For a set $\cP$ of phylogenetic trees, we define
the {\em excess of $\cP$}, denoted $\exc(\cP)$ by
$$\exc(\cP)=|L(\cP)|-3-\sum_{T \in \cP}|\Er(T)|.$$ 
Finally, we say that $\cP$ is {\em excess-free} if $\exc(\cP)=0$, we define 
$\Exc(\cP)=\{\cP' \subseteq \cP:\exc(\cP')=0\}$, and we call $\cP$ {\em slim} if $\exc(\cP') \ge 0$ for
every non-empty subset $\cP' \subseteq \cP$. In the next section, we
will prove: 
\begin{theorem}
Every slim set of binary phylogenetic trees is compatible. 
\label{slim}
\end{theorem}
It depends only on the involved leaf sets of the trees and not on
which phylogenetic tree is chosen for a fixed leaf set whether a
collection of binary phylogenetic trees is slim. Nevertheless, Theorem
\ref{slim} is sharp in the sense that, for every binary
phylogenetic tree $T$, there are slim sets of quartets that define
$T$ \cite{bds99}. 

After proving Theorem \ref{slim} in Section~\ref{proof}, we will apply
it in Section~\ref{appl} to give a
new proof of the known characterization of excess-free definitive
quartet sets. An outline of the original proof is given in
\cite{bds99} while the full details can be found in \cite{b99}. 
That characterization has been referred to as ``one of the most
mysterious and apparently difficult results in phylogenetics''
\cite{ds04}. 

\section{Proof of Theorem \ref{slim}}
\label{proof}

We assume that the statement is wrong and that $\cP$ is a minimal
non-compatible slim set of binary phylogenetic trees where $L(\cP)$ is
minimal under
all non-compatible slim sets of binary phylogenetic trees with cardinality
$|\cP|$. Clearly, we have $|L(T)| \ge 4$ for every tree $T \in \cP$ and
there is no subset $\cP_1$ of $\cP$ with $\exc(\cP_1)=0$ and
$1<|\cP_1|<|\cP|$ since otherwise $\cP_1$ can be replaced by a single
binary phylogenetic tree with leaf set $L(\cP_1)$ that displays
$\cP_1$, contradicting the minimality of $|\cP|$. 

We claim that, for every  cherry $\{x,y\}$ of $T$ for some tree $T$ in
$\cP$, there is another tree $T'$ in $\cP$ that contains $x$ and $y$ as
leaves but not as a cherry. Assuming that the claim is wrong, let
$\cP'$ be the set of phylogenetic trees obtained from $\cP$ by
replacing $x$ and
$y$ by a new leaf $z$, i.e. if a tree $T \in \cP$ contains either $x$
or $y$ then this leaf is replaced by $z$ and if $T$ contains a cherry
$\{x,y\}$ then the cherry is replaced by $z$. Clearly, we have
$|L(\cP')|<|L(\cP)|$ and $\cP'$ is not 
compatible since the tree obtained from a phylogenetic tree displaying
$\cP'$ by replacing $z$ by a cherry $\{x,y\}$ would display $\cP$. Hence,
by induction hypothesis, $\cP'$ is not slim. Let $\cP_2'$ be a subset of
$\cP'$ with negative excess and let $\cP_2$ be the subset of $\cP$
obtained from $\cP_2'$ by reversing the identification of $x$ and $y$,
i.e. replacing $z$ by a cherry $\{x,y\}$ or a
leaf $x$ or $y$ for
every tree in $\cP_2'$ that contains $z$. We have 
\begin{eqnarray*}\exc(\cP_2) & = & |L(\cP_2)|-3-\sum_{T \in
    \cP_2}|\Er(T)|  \\
& \le & |L(\cP_2')|+1-3-\sum_{T \in \cP_2'}|\Er(T)| \\ 
& = & \exc(\cP_2')+1 \le 0.
\end{eqnarray*}
This implies $\cP_2=\cP$ or $|\cP_2|=1$ but for both cases we get 
$\exc(\cP_2') \ge \exc(\cP_2) \ge 0$, a contradiction. This proves the claim. 

We construct a digraph $G_{\cP}=(V_{\cP},A_{\cP})$ with vertex set
$V_{\cP}=\cP$ and an arc labeling $\lambda$ that associates a
set of two elements of $L(\cP)$ with every arc $a \in A_{\cP}$ as follows: 
For every tree $T \in \cP$, we choose two different cherries
$\{x_1(T),y_1(T)\}$ and $\{x_2(T),y_2(T)\}$ of $T$ and, for $i \in \{1,2\}$,
we choose $T_i(T) \in \cP$ such that $x_i(T),y_i(T) \in L(T_i)$ but
$\{x_i(T),y_i(T)\}$ is not a cherry of $T_i(T)$. Then we define 
$A_{\cP}=\{(T,T_i(T)):T \in \cP, i \in \{1,2\}\}$ and
$\lambda((T,T_i(T)))=\{x_i(T),y_i(T)\}$ for all $T \in \cP$ and $i \in
\{1,2\}$. For the ease of notation, we will write $\lambda(T_1,T_2)$
rather than $\lambda((T_1,T_2))$ for an arc $(T_1,T_2) \in A_{\cP}$. 

\begin{figure}[h]
\begin{center}
\input{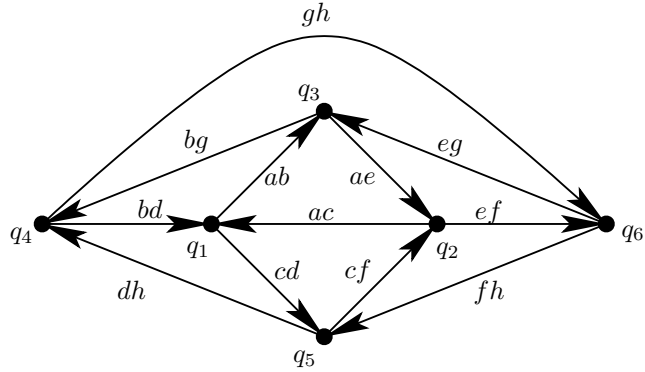}
\caption{The (in this case well-defined) arc-labeled digraph $G_{\cP}$ for $\cP=\{q_1,...,q_6\}$ with $q_1=ab|cd$, $q_2=ac|ef$, $q_3=ae|bg$, $q_4=bd|gh$, $q_5=cf|dh$, and $q_6=eg|fh$.}
\label{octa}
\end{center}
\end{figure}

By construction, the outdegree $d^+(T)=2$ for every tree $T \in
\cP$ thus we have $|A(G')|\le 2|V(G')|$ for every subdigraph $G'$ of
$G_{\cP}$ and equality holds for $G'=G_{\cP}$. Let $G$ be a minimal
subdigraph of $G_{\cP}$ with $|A(G)|=2|V(G)|$, let $\cP_3$ be the
vertex set of $G$
and let $X=L(\cP_3)$. We define an {\em
  $x$-coloured cycle} of $G$ to be a cycle $C$ of $G$ such that $x \in
\lambda(a)$ for every arc $a$ of $C$. We call $C$ a {\em coloured
  cycle} if there is $x \in X$ such that $C$ is an $x$-coloured cycle. 

Since two different cherries of a binary tree are disjoint,
$\lambda(T_1,T_2) \cap \lambda(T_1,T_3) \neq \emptyset$ implies
$T_2=T_3$. Hence, every coloured cycle $C$ of $G$ is a directed cycle
since otherwise there is a vertex of outdegree 2 in $C$ but the
two outgoing arcs must have disjoint label sets. This implies that every
arc with label set $\{x,y\}$ is contained in at most one $x$-coloured
and in at most one $y$-coloured cycle. 

Recall that the {\em cycle space} of a graph $H=(V,E)$ is the vector space over the 2-element field 
${\mathbf Z}_2$ generated by the cycles of $H$ (considered as binary vectors indexed by $E$), together with the symmetric difference as addition. This makes it possible to decide whether a collection of cycles of $H$ is linearly dependent. The dimension of the cycle space is called the {\em cyclomatic number} 
$c(H)$, and if $H$ is connected we have $c(H)=|E|-|V|+1$. For a digraph, we define the cycle space and the cyclomatic number to be the ones of the underlying undirected graph. 

We will now show that the number of coloured cycles in $G$ is at most
$c(G)+1$. This 
clearly holds if all coloured cycles of $G$ are linearly independent. 
Let $S$ be a minimal linearly dependent set of coloured cycles of $G$
and let $G_S$ be the subdigraph of $G$ induced by $V(G_S)= \bigcup_{C
  \in S} V(C)$. By the minimality of $S$, $G_S$ is connected and every
arc of $G_S$ is contained in at least two
cycles of $S$. Let $(T_1,T_2)$ and $(T_2,T_3)$ be arcs of an
$x$-coloured cycle $C \in S$, let $\lambda(T_1,T_2)=\{x,y\}$ and
$\lambda(T_2,T_3)=\{x,z\}$. We have $y \neq z$ and there are arcs
$(T_4,T_2)$ and $(T_2,T_5)$ in $G_S$ with $y \in \lambda(T_2,T_5)$ and
$z \in \lambda(T_4,T_2)$ since there is a $y$-coloured cycle
containing $(T_1,T_2)$ and a $z$-coloured cycle containing $T_2,T_3$
in $S$. Hence, the degree of $T_2$ in $G_S$ as well as the degree of
every other vertex of $G_S$ is at least 4. This implies $|A(G_S)| \ge
2|V(G_S)|$ but we have $|A(G')| < 2|V(G')|$ for every proper
subdigraph of $G$. So we must have $G_S=G$ and $S$ is the set of all
coloured cycles of $G$. This means that the number of coloured cycles
in $G$ is at most $c(G)+1$. Since $G$ is connected we have
$c(G)=|A(G)|-|V(G)|+1=2|\cP_3|-|\cP_3|+1=|\cP_3|+1$ and the number of
coloured cycles in $G$ is at most $|\cP_3|+2$. 

We conclude the proof by showing that $\cP_3$ has negative excess, in contradiction to our assumption that $Q$ is slim. For
$x \in X$, let $G_x=(V_x,A_x)$ be the subdigraph of $G$ such that $V_x$ is the set of all trees in 
$\cP_3$ that contain $x$ and $A_x$ is the set of all arcs of $G$ whose label sets contain $x$. Then 
all cycles in $G_x$ are vertex-disjoint, so we have $|V_x| \ge 1 + |A_x| - c_x$ where $c_x$ is the number of cycles in $G_x$. 
This implies $$\sum_{x \in X}|V_x| \ge |X|+2|A(G)|-|\cP_3|-2=|X|+3|\cP_3|-2.$$
Finally, we get 
\begin{eqnarray*}
\exc(\cP_3) & = & |X|-3- \sum_{T \in \cP_3}|\Er(T)| \\
& = & |X|-3-\sum_{T \in \cP_3}(|L(T)|-3) \\
& = & |X|-3+3|\cP_3|-\sum_{x \in X}|V_x| \\
& \le & |X|-3+3|\cP_3|-|X|-3|\cP_3|+2 \\ 
& = & -1. 
\end{eqnarray*} \epf

\section{Definitive Quartet Sets of Minimum Size}
\label{appl}
Let $T$ be a phylogenetic tree with
$\{a,b,c,d\} \subseteq
L(T)$ and let $e \in \Er$. Then the quartet $ab|cd$ {\em
  distinguishes} $e$ if $T$ displays $ab|cd$ and $e$ is the only edge that the paths from
$a$ to $c$ and from $b$ to $d$ in $T$ have in common. 
We summarize some elementary results
about definitive quartet sets. Proofs can be found in \cite{ss03}. 
\begin{lemma}
\label{easy}
Let $Q$ be a quartet set that defines a phylogenetic tree $T$ with
$|L(T)|=n$. Then $T$ is binary and for every interior edge
$e$ of $T$, there must be a quartet in $Q$ that distinguishes $e$. Hence, we
have $|Q| \ge n-3$ and, in case of equality, every quartet in $Q$ distinguishes an 
interior edge of $T$. 
\end{lemma}

Combining Theorem \ref{slim} with Lemma \ref{easy}, we get the following corollary: 

\begin{corollary}
\label{definitive}
Let $Q$ be a quartet set that defines a phylogenetic tree $T$ with $|L(T)|=n$. If $|Q|=n-3$, then $Q$ is slim and every excess-free subset of $Q$ is definitive. 
\end{corollary}
\pf
Assume that $|Q|=n-3$ and $Q$ is not slim. Let $Q'$ be a subset of $Q$ with 
negative excess. Since $T|_{L(Q')}$ contains $|L(Q')|-3$ 
interior edges but $Q'$ contains more than $|L(Q')|-3$ quartets there 
must be an interior edge in $T$ that is not distinguished by a quartet 
in $Q$, a contradiction to Lemma \ref{easy}. Now assume that an excess-free subset $Q'$ of $Q$ is not definitive. Let $T'$ with $L(T')=L(Q')$ be a different tree from $T|_{L(Q')}$ that displays $Q'$. Then $Q-Q'+T'$ is slim and compatible in view of Theorem \ref{slim}. However, a tree that displays $Q-Q'+T'$ displays $Q$ as well but it is different from $T$. Hence, $Q$ can not be definitive.  
\epf

For a set of two binary phylogenetic trees, it is easy to decide if
they are compatible resp. definitive. The proof of the next lemma is
straight forward and omitted. 

\begin{lemma}
\label{2trees}
Let $\cP=\{T_1,T_2\}$ be a set of two binary phylogenetic trees with
$L(T_1)\cap L(T_2) \neq \emptyset$. Then $\cP$ is compatible if and
only if there is a binary tree $T_{1,2}$ with leaf set $L(T_1) \cap
L(T_2)$ that is displayed by $T_1$ and $T_2$. $\cP$ is
definitive if and only if, in addition, the sets of edges in
$T_{1,2}$ which are subdivided in $T_1$ and $T_2$, respectively, are
disjoint. 
\end{lemma}

\begin{figure}[h]
\begin{center}
\input{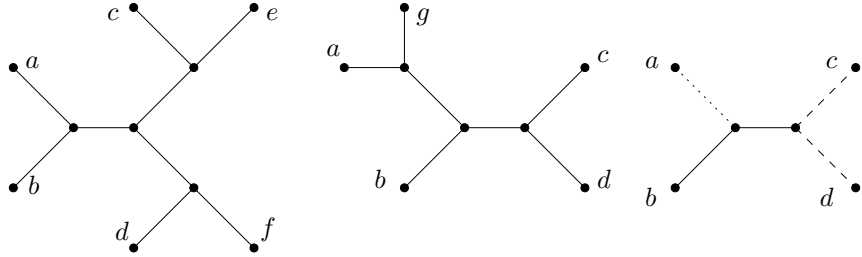}
\caption{Two trees $T_1$ and $T_2$ and the restriction $T_{1,2}$ of both trees to the intersection of their leaf sets. The dashed edges have to be subdivided to obtain $T_1$ from $T_{1,2}$ and the dotted edge has to be subdivided to obtain $T_2$.}
\label{component}
\end{center}
\end{figure}

For a set $Q$ and a set $\cC$ of subsets of $Q$, the
pair $(Q,\cC)$ is called a {\em patchwork} if, for every two sets
$C_1,C_2 \in \cC$ with $C_1 \cap C_2 \neq \emptyset$, we have $C_1\cap
C_2 \in \cC$ and $C_1 \cup C_2 \in \cC$. Patchworks were introduced
by B\"ocker and Dress \cite{bd01}. An easy
proof of the following lemma can be found in \cite{ss03}, p.~136.  

\begin{lemma}
Let $T$ be a binary phylogenetic tree and let $Q$ be an excess-free
quartet set such that every interior edge of $T$ is distinguished
by an element of $Q$. Then $(Q,\Exc(Q))$ is a patchwork. 
\label{patchwork}
\end{lemma}

A set $\cC$ of non-empty subsets of a set $Q$ is called a
{\em hierarchy} of $Q$ if, for every two sets $C_1,C_2 \in \cC$ with
$C_1 \cap C_2 \neq \emptyset$, we have $C_1 \subseteq C_2$ or $C_2
\subseteq C_1$. The elements of a hierarchy are called {\em
  clusters}. A hierarchy $\cC$ of $Q$ is {\em maximal} if $Q \in
\cC$ and, for every set $C \in \cC$ with $|C|>1$, there are $C_1,C_2
\in \cC$ such that $C=C_1 \cup C_2$ and $C_1 \cap C_2= \emptyset$. We
call $C_1$ and $C_2$ the
{\em children} of $C$. This is only one of several different
definitions of a maximal hierarchy and their equivalence is shown in
\cite{bd00}. A patchwork $(Q,\cC)$ is called
{\em ample} if there is a maximal hierarchy $\cC'$ of $Q$ with $\cC'
\subseteq \cC$. The next lemma is a special case of the results in
\cite{bd01} where ample patchworks are studied in detail. 

\begin{lemma}
Let $Q$ be a finite set and let $(Q,\cC)$ be an ample patchwork with
$C \in \cC$. Then there is a maximal hierarchy $\cC_C$ of $Q$ with $C
\in \cC_C \subseteq \cC$. 
\label{ample}
\end{lemma}
\pf 
Assume the statement is wrong and $C \in \cC$ is maximal such that there is no maximal hierarchy $\cC_C$ of $Q$ with $C \in \cC_C \subseteq \cC$. Clearly, we have $C \neq Q$. Let $\cC' \subseteq \cC$ be a maximal hierarchy of
$Q$. Then $\cC_1=\{C \cap C':C' \in \cC' \mbox{ and } C \cap C' \neq
\emptyset\}$ is a maximal hierarchy of $C$ with $\cC_1 \subseteq
\cC$.  Let $A \in \cC'$ be minimal such that $A \cap C \neq \emptyset$ and $A \cup C \neq C$ and let $B$ be the child of $A$ with $B \cap C = \emptyset$. Since $\cC$ is a patchwork we have that $A \cup C \in \cC$ and, by definition, $A \cup C$ is the disjoint union of $B$ and $C$. In view of the maximality of $C$ there is a maximal hierarchy $\cC'' \subseteq \cC$ of $Q$ with $A \cup C \in \cC''$. By replacing all proper subsets of $A \cup C$ in $\cC''$ by the union of $\cC_1$ and the maximal hierarchy of $B$ contained in $\cC'$, we obtain a maximal hierarchy $\cC_C$ of $Q$ with $C \in \cC_C \subseteq \cC$, in contradiction to our assumption. 
\epf

Combining Lemma \ref{patchwork} and Lemma \ref{ample} yields the
following result. 

\begin{corollary}
\label{amplecor}
Let $Q$ be a set of quartets with $C \in \Exc(Q)$. If $\Exc(Q)$ contains
a maximal hierarchy of $Q$, then $\Exc(Q)$ contains a maximal hierarchy
$\cC_C$ of $Q$ with $C \in \cC_C$. 
\end{corollary}

We are now able to establish the hardest implication of the
characterization of excess-free definitive quartet sets. 

\begin{theorem} 
\label{hard}
If an excess-free set $Q$ of quartets defines a binary phylogenetic
tree $T$, then $\Exc(Q)$ contains a maximal hierarchy of $Q$. 
\end{theorem} 
\pf
To simplify the notation, we call a maximal hierarchy of a subset $Q_s
\subseteq Q$ that is contained in $\Exc(Q)$ a {\em $Q_s$-hierarchy}. 

Assume that $Q$ is a minimal counterexample where $\{w_1,x_1\}$ is a cherry
of $T$ and $w_1x_1|w_2y \in Q$. Note that $Q$ is slim in view of Corollary \ref{definitive}. 
Then there is a subset $Q'$ of
$Q-w_1x_1|w_2y$ such that $Q'$ contains $w_2$ and $y$ and $\exc(Q')=0$
since otherwise, by Theorem \ref{slim}, the quartet set obtained from
$Q$ by identifying $w_2$ and $y$ would be compatible and thus there
would be a phylogenetic tree displaying $Q$ that contains a cherry
$\{w_2,y\}$. However, $\{w_2,y\}$ is not a cherry of $T$ since $w_1x_1|w_2y$ distinguishes 
an edge of $T$. 

For every subset $Q_s$ of $Q-w_1x_1|w_2y$, we let 
$Q_s^*$ denote the quartet set obtained from 
$Q_s$ by replacing $w_1$ and $x_1$ by a new leaf
$z_1$ and we define $Q^*=(Q-w_1x_1|w_2y)^*$. Conversely, for $Q_{s'}^*\subseteq Q^*$, 
we use $Q_{s'}$ to denote the quartet set obtained from $Q_{s'}^*$ by reversing the replacement. Note that $Q_{s'}$ is well defined since $Q$ is excess-free and defines $T$. Hence, every quartet in $Q$ distinguishes an interior edge of $T$ and every interior edge of $T$ is distinguished by exactly one quartet. 
Clearly, $Q^*$ defines the tree $T^*$ obtained from $T$ by replacing the cherry $\{w_1,x_1\}$ by $z_1$. Further, since we have $|Q^*|\le |L(T^*)|-3$ and Lemma \ref{easy} implies $|Q^*| \ge |L(T^*)|-3$, 
$Q^*$ is excess-free thus, by induction hypothesis,
$(Q^*,\Exc(Q^*))$ is an ample patchwork. 

Since $\exc(Q')=0$ and $Q$ is slim we have $\exc(Q'^*)=0$ thus, by
Corollary \ref{amplecor}, there is a $Q^*$-hierarchy that contains $Q'^*$. Hence, one child
$Q_1^*$ of $Q^*$ contains $w_2$ and $y$. Let $Q_2^*$ be the other
child of $Q^*$. We note that $|L(Q_1^*) \cap L(Q_2^*)|=3$.  

We will now consider the case $z_1 \in L(Q_1^*) \cap L(Q_2^*)$. Let $L(Q_1^*) \cap L(Q_2^*)=\{z_1,z_2,z_3\}$. Then $Q_1+w_1x_1|w_2y$ is excess-free and in view of Corollary \ref{definitive} definitive. By induction hypothesis, there is a $(Q_1+w_1x_1|w_2y)$-hierarchy $\cC_1$. Since $Q$ is definitive
$Q_2+w_1x_1|z_2z_3$ is definitive because every tree that displays $Q_2+w_1x_1|z_2z_3$ and the tree defined by $Q_1+w_1x_1|w_2y$ are compatible in view of Lemma \ref{2trees}. Therefore,  there is a $(Q_2+w_1x_1|z_2z_3)$-hierarchy $\cC_2$. Let $\cC_2'$ be the set system that is obtained from $\cC_2$ by replacing every cluster $C$ that contains $w_1x_1|z_2z_3$ by $(C-w_1x_1|z_2z_3) \cup (Q_1+w_1x_1|w_2y)$. Then $\cC_1 \cup \cC_2'$ is a $Q$-hierarchy, contradicting our assumption. Hence, we have $z_1 \notin L(Q_1^*) \cap L(Q_2^*)$.

If $z_1 \in L(Q_1^*)$  and $z_1 \notin L(Q_2^*)$ then $\exc(Q_1+w_1x_1|w_2y)=\exc(Q_2)=0$ thus there is a $Q$-hierarchy where $Q_1+w_1x_1|w_2y$
and $Q_2=Q_2^*$ are the children of $Q$. Hence, we have $z_1 \in L(Q_2^*)$  and $z_1 \notin L(Q_1^*)$. We can assume that $w_2 \notin
L(Q_2^*)$ since in case $w_2,y \in L(Q_2^*)$ interchanging $Q_1^*$ and
$Q_2^*$ would yield $\{w_2,y,z_1\} \subset L(Q_1^*)$, a case that we have already considered. If $y \notin
L(Q_2^*)$ then let $L(Q_1) \cap L(Q_2)=\{v_1,v_2,v_3\}$. Lemma \ref{2trees} implies that different edges of the phylogenetic tree with leaf set $\{v_1,v_2,v_3\}$ are subdivided in $T|_{L(Q_1)}$ and $T|_{L(Q_2)}$ but then $w_1x_1|w_2y$ does not distinguish an edge of $T$, a contradiction to Lemma \ref{easy}. 

The only remaining case is $w_2,y \in L(Q_1^*)$, $z_1 \notin
L(Q_1^*)$, $y,z_1 \in L(Q_2^*)$, and $w_2 \notin
L(Q_2^*)$. Let $L(Q_1^*) \cap L(Q_2^*) =
\{y,x_2,w_3\}$. In view of Lemma \ref{2trees} the edge sets in the phylogenetic tree with leaf set 
$\{y,x_2,w_3\}$ that have to be subdivided to obtain $T|_{L(Q_1)}$ respectively $T|_{L(Q_2)}$ are disjoint. On the other hand, $\{w_1,x_1\}$ is a cherry of $T|_{L(Q_2)}$ while $w_2$ is in $L(Q_1)$. Since 
$w_1x_1|w_2y$ distinguishes an edge of $T$ it follows that 
$w_2y|x_2w_3$ is not displayed by $T$. Hence, we can assume without loss of generality that $T$ displays $w_2x_2|w_3y$. 
Further, we have $Q_1=Q_1^*$ thus there is a $Q_1$-hierarchy. 
Let $T_2^*$ be the tree defined by $Q_2^*$. Then the edge of the phylogenetic tree with leaf set 
$\{x_2,w_3,y\}$ that is incident with $w_2$ does not have to be subdivided to obtain $T_2^*$. 
Therefore, $Q_2^*+w_2x_2|yw_3$ is definitive. Further, the quartet set 
$Q_2^+:=Q_2+w_2x_2|yw_3+w_1x_1|w_2y$ is definitive since every tree that displays $Q_2^+$ 
is compatible with the tree defined by $Q_1$ in view of Lemma \ref{2trees}, and $Q_2^+$ is excess-free. 
If  $|Q_1|>1$, then $|Q_2^+|<Q$ thus, by induction hypothesis, there is a
$(Q_2+w_1x_1|w_2y+w_2x_2|w_3y)$-hierarchy and we can construct a $Q$-hierarchy in the same way as for the case $z_1 \in L(Q_1^*) \cap L(Q_2^*)$. Hence,
$Q_1=\{w_2x_2|w_3y\}$ thus $\{w_2,x_2\}$ is a cherry of $T$ and
$w_1x_1|w_2y$ and $w_2x_2|w_3y$ are the only two quartets in $Q$ which
contain $w_2$. 

We have shown that, for every cherry $\{w,x\}$ of $T$ with $wx|w'y'
\in Q$, there is $z' \in \{w',y'\}$ such that the following statements hold:  
\begin{enumerate}
\item [(P1)]
$z'$ is contained in a cherry $\{z',x'\}$ of $T$.  
\item[(P2)]
There are exactly two quartets in $Q$ which contain $z'$.  
\item [(P3)]
There is $w'' \in L(T)$ such that $w'x'|w''y' \in Q$ (if $z'=w'$) or $x'y'|w'w'' \in Q$ (if $z'=y'$). 
\end{enumerate}
We will use Properties (P1)--(P3) of $Q$ and $T$ to show that, for every natural number $i$, the following statements hold: 
\begin{enumerate}
\item [(P4)]
There is a leaf set $L_i=\{w_j:1 \le j \le i+1\}\cup \{x_j:1 \le j \le i\} \cup \{y\} \subset L(T)$ of cardinality $2i+2$ such that $T|_{L_i}$ is the phylogenetic tree $T_i$ with vertex set $L_i \cup \{u_j,v_j:1 \le j \le i\}$ and edge set $\{ yu_1\}\cup\{u_jv_j,v_jw_j,v_jx_j,u_ju_{j+1}:1\le j\le i\}$ where we define $u_{i+1}=w_{i+1}$. 
\item [(P5)]
$\{w_i,x_i\}$ is a cherry of $T$. 
\item[(P6)]
$w_ix_i|w_{i+1}y \in Q$.
\end{enumerate}
This will complete the proof of the theorem since it contradicts the finiteness of $T$. 

We have already shown that all claimed statements hold for $i \le 2$. Assume $i \ge 3$ and that all claimed statements are true for smaller $i$. Since, by induction hypothesis, $\{w_{i-1},x_{i-1}\}$ is a cherry of $T$ in view of (P5)  and $w_{i-1}x_{i-1}|w_iy \in Q$ in view of (P6), (P1) implies that there is a leaf $x_i$ of $T$ and $z \in \{w_i,y\}$ such that $\{z,x_i\}$ is a cherry of $T$. Since $T|_{L_{i-1}}=T_{i-1}$ we have $x_i \notin L_{i-1}$. By (P2), there are exactly two quartets in $Q$ which contain $z$. Since $y$ is contained in the quartets $w_1x_1|w_2y$ and $w_2x_2|w_3y$ and there is a quartet containing $x_i$ and $y$ if $\{x_i,y\}$ is a cherry of $T$ we can not have $z=y$. 
Hence, $z=w_i$ and (P3) implies that there is $w_{i+1} \in L(T)$ with $w_ix_i|w_{i+1}y \in Q$.  Since $w_ix_i|w_{i+1}y$ distinguishes an edge of $T$ we have $T|_{L_i}=T_i$. 
\epf

The characterization of excess-free defining quartet sets which was originally published in \cite{bds99} is an easy
consequence of Theorem~\ref{hard}. 

\begin{corollary}
\label{char}
Let $T$ be a binary phylogenetic $X$-tree and let $Q$ be an
excess-free quartet set displayed by $T$ with $L(Q)=L(T)$. Then $Q$
defines $T$ if and only if (i)~every interior edge of $T$ is
distinguished by a quartet in $Q$ and (ii)~$\Exc(Q)$ contains a
maximal hierarchy. 
\end{corollary}
\pf
If $Q$ defines $T$ then (i) follows from Lemma \ref{easy} and (ii)
from Theorem \ref{hard}. Conversely, we use induction on $|L(T)|$ to
show that $Q$ defines $T$ if (i) and (ii) hold. This is clearly true
if $|L(Q)|=4$ thus we can assume that $|L(Q)| \ge 5$ and there is a
maximal hierarchy contained in $\Exc(Q)$ where $Q_1$ and $Q_2$
are the children of $Q$. By induction hypothesis, $Q_1$ and $Q_2$
define binary trees $T_1$ and $T_2$, respectively. For $i \in
\{1,2\}$, every quartet in $Q_i$ distinguishes an interior edge of
$T$. Hence, the interior edges of $T_i$ are not subdivided in the
smallest subtree of $T$ containing
$L(T_i)$. Therefore, the interior edges of $T$ distinguished by a
quartet in $Q_i$ induce a connected
subgraph $T_i'$ of $T$ which implies that $T_1'$ and $T_2'$ have
exactly one vertex $v$ in common. Let $L(T_1) \cap
L(T_2)=\{x,y,z\}$. By symmetry, we can assume that $\{x,y\}$ is a
cherry of $T_1$ and $z$ is an isolated vertex in the graph obtained
from $T_2$ by removing the path from $x$ to $y$. Hence, only the edge of the phylogenetic tree with leaf set $\{x,y,z\}$ that is incident with $z$ has to be subdivided to obtain $T_1$ while one or both of the other edges have to be subdivided to obtain $T_2$. Therefore, Lemma \ref{2trees} implies that $T_1$ and $T_2$ define $T$, thus $Q$ defines $T$. 
\epf

We conclude the paper by pointing out one of the algorithmic
consequences of Corollary \ref{char} that are studied in
\cite{bbds00}. 

\begin{corollary} 
It can be decided in polynomial time if an excess-free quartet set $Q$
is definitive and then the phylogenetic tree displaying $Q$ can be
reconstructed in polynomial time. 
\end{corollary}

\section*{Acknowledgment} 
I thank the  Allan Wilson Centre for Molecular Ecology and Evolution and the Department of Mathematics and Statistics, University of Canterbury in Christchurch, New Zealand where I worked when I found the presented results and I thank Mike Steel for fruitful discussions. I would like to thank Meng Chen and Katharina Huber for their comments on earlier versions of this paper.

\end{document}